\documentclass[reqno, 11pt]{article}

\usepackage{enumerate}
\usepackage{needspace}
\usepackage{latexsym}
\usepackage{amsmath}
\usepackage{amsfonts}
\usepackage{amsthm}
\usepackage{amssymb,mathtools}
\usepackage{newlfont}
\usepackage{graphics}
\usepackage{color}
\usepackage{float}
\usepackage{diagbox}
\usepackage{tocloft}
\usepackage{titlesec}
\usepackage{booktabs,array,longtable,enumitem,microtype}
\usepackage{extpfeil}
\usepackage{centernot}

\usepackage[linesnumbered,ruled,vlined]{algorithm2e}
\usepackage{url}

\usepackage[margin=1in]{geometry}

\usepackage[T1]{fontenc}
\usepackage{lmodern}
\usepackage{microtype}
\usepackage{rotating}
\usepackage{multirow}
\usepackage{extarrows}
\usepackage[sort,compress,numbers]{natbib}
\usepackage[utf8]{inputenc}
\usepackage{xcolor}
\usepackage{listings}
\usepackage{aliascnt}
\usepackage[pagebackref,colorlinks=true,linkcolor=blue,citecolor=red,urlcolor=blue]{hyperref}
\usepackage[nameinlink,noabbrev,capitalize]{cleveref}
\numberwithin{equation}{section}

\setlist[enumerate]{label=\textup{(\roman*)},leftmargin=*}
\newtheorem{theorem}{Theorem}[section]
\newaliascnt{proposition}{theorem}
\newtheorem{proposition}[proposition]{Proposition}
\aliascntresetthe{proposition}
\newaliascnt{lemma}{theorem}
\newtheorem{lemma}[lemma]{Lemma}
\aliascntresetthe{lemma}
\newaliascnt{corollary}{theorem}
\newtheorem{corollary}[corollary]{Corollary}
\aliascntresetthe{corollary}

\theoremstyle{definition}
\newaliascnt{definition}{theorem}

\aliascntresetthe{definition}
\newaliascnt{remark}{theorem}
\newtheorem{remark}[remark]{Remark}
\aliascntresetthe{remark}
\newaliascnt{example}{theorem}

\aliascntresetthe{example}
\newaliascnt{conjecture}{theorem}
\newtheorem{conjecture}[conjecture]{Conjecture}
\aliascntresetthe{conjecture}
\crefname{conjecture}{Conjecture}{Conjectures}
\newaliascnt{question}{theorem}

\aliascntresetthe{question}
\newaliascnt{problem}{theorem}

\aliascntresetthe{problem}
\newaliascnt{alg}{theorem}

\aliascntresetthe{alg}

\allowdisplaybreaks[4]

\SetKwInput{KwInput}{Input}                
\SetKwInput{KwOutput}{Output}              

\newcommand{\R}{\mathbb R}
\newcommand{\Z}{\mathbb Z}
\DeclareMathOperator{\Int}{int}
\DeclareMathOperator{\vol}{vol}
\DeclareMathOperator{\area}{area}
\DeclareMathOperator{\lw}{lw}
\DeclareMathOperator{\conv}{conv}
\numberwithin{equation}{section}

\title{The Zaks--Perles--Wills volume conjecture I:\\ three-dimensional lattice polytopes}

\author{Feihu Liu$^{\color{blue} \dag}$, Jinlong Tang$^{\color{blue} \ddag}$, Sihao Tao$^{\color{blue} \P}$, and Zihao Zhang$^{\color{blue} \S}$
\\[2mm]
{\small $^{\color{blue} \dag}$ Center for Combinatorics, LPMC,}\\[-0.8ex]
{\small Nankai University, Tianjin 300071, P.R.~China}\\
{\small $^{\color{blue} \ddag, \P}$ School of Mathematical Sciences,}\\[-0.8ex]
{\small Capital Normal University, Beijing, 100048, P.R.~China}\\
{\small $^{\color{blue} \S}$ School of Mathematics and Statistics,}\\[-0.8ex]
{\small Beijing Institute of Technology, Beijing 102400, P.R.~China}\\
{\small {\color{blue} $^\dag$} Email address: \url{liufeihu7476@163.com}}\\
{\small {\color{blue} $^\ddag$} Email address: \url{jinlong\_tang@cnu.edu.cn}}\\
{\small {\color{blue} $^\P$} Email address: \url{sihao\_tao@cnu.edu.cn}}\\
{\small {\color{blue} $^\S$} Email address: \url{zihao-zhang@foxmail.com}}\\
}

\date{\today}

\begin{document}

\maketitle

\begin{abstract}
We prove the three-dimensional Zaks--Perles--Wills volume bound for lattice polytopes with $k\ge3$ interior lattice points: the Euclidean volume of such a lattice polytope is at most $6(k+1)$, with equality if and only if the polytope is affinely unimodularly equivalent to $\operatorname{conv}(0,2e_1,3e_2,6(k+1)e_3)$. The inequality is strict when the minimum lattice width is at least three. Together with the known cases $k=1,2$, this completely resolves the three-dimensional Zaks--Perles--Wills volume conjecture.
\end{abstract}

\noindent
\begin{small}
\end{small}

\noindent
\begin{small}
\end{small}


\section{Introduction}\label{sec:introduction}

A \emph{lattice polytope} in $\R^d$ is the convex hull of finitely many points of $\Z^d$.
Unless stated otherwise, lattice polytopes are assumed to be full-dimensional.
For such a polytope $P$, write
\[
I(P)=|\Int(P)\cap\Z^d|,\qquad
\operatorname{Vol}(P)=d!\vol(P),
\]
where $\Int(P)$ is its interior, $\vol(P)$ is its Euclidean volume, and
$\operatorname{Vol}(P)$ is its normalized volume.
Two polytopes $P,Q\subset\R^d$ are \emph{affinely unimodularly equivalent},
written $P\simeq Q$, if $Q=UP+b$ for some $U\in\mathrm{GL}_d(\Z)$
and $b\in\Z^d$. Such transformations preserve the lattice, volume,
and the number of interior lattice points.

The volume of lattice polytopes has been studied from both computational and geometric perspectives; for computational aspects, see
\cite{BBCVVolum,DeLoera}. A related problem is Ehrhart's volume conjecture
\cite{EhrhartVolume,BNillPaffenhonlz}: if a full-dimensional compact convex
body $K\subset\R^d$ has barycenter $0$ and
$\Int(K)\cap\Z^d=\{0\}$, then
\[
\vol(K)\le\frac{(d+1)^d}{d!},
\]
with equality precisely when  $
 K=U\bigl((d+1)\conv(0,e_1,\ldots,e_d)-(1,\ldots,1)^{\mathsf{T}}\bigr)$ for some $U\in\operatorname{GL}_d(\mathbb Z)$,
where $e_1,\ldots,e_d$ are the standard basis vectors of $\mathbb R^d$.
The volume bound was proved by OpenAI \cite[Chapter~8, Theorem~1.1]{OAI26},
and the equality case was established by Liu \cite[Theorem~1.3]{LiuJiHao}.
For further volume bounds, see \cite{Bau25,Nill07}.

For fixed dimension and a fixed positive number of interior lattice points,
Hensley \cite[Theorem~3.6]{Hensley1983} proved that the volume of a lattice
polytope is bounded. Quantitative estimates were subsequently obtained by
Pikhurko \cite{Pikhurko2001}. The expected sharp bound is attained by the
simplices constructed by Zaks, Perles, and Wills \cite{ZPW1982}.

To define these simplices, let $s_0=1$, $s_1=2$ and
$s_{j+1}=s_1\cdots s_j+1$ for $j\ge1$. These are the terms of the
\emph{Sylvester sequence} \cite[OEIS; A000058]{Sloane23}.  Thus we have
$$(s_n)_{n\geq 0} = (1, 2, 3, 7, 43, 1807, 3263443, \ldots)$$
For $d\ge2$ and $k\ge1$, Zaks, Perles, and Wills ~\cite{ZPW1982} introduced the $d$-dimensional simplex:
\[
S_{d,k}=\conv\bigl(0,s_1e_1,\ldots,s_{d-1}e_{d-1},(k+1)(s_d-1)e_d\bigr).
\]
These simplices are called \emph{Sylvester simplices} \cite{CCJP2026}, or \emph{Zaks--Perles--Wills simplices} \cite{Ave20}. 
The Ehrhart properties of Sylvester simplices have been studied by Caicedo, Castillo, Juhnke, and Poullot in \cite{CCJP2026}.
The interior lattice points of $S_{d,k}$ are exactly
$(1,\ldots,1,j)$ for $1\le j\le k$. This follows directly from
the defining inequalities and the identity
$\sum_{i=1}^{d-1}\frac{1}{s_i}=1-\frac{1}{s_d-1}$. Thus
\[
I(S_{d,k})=k,\qquad
\vol(S_{d,k})=\frac{(k+1)(s_d-1)^2}{d!}.
\]

The following conjecture goes back to Zaks, Perles, and Wills \cite{ZPW1982}; the explicit formulation below appears in
\cite[Conjecture~1.5]{BK2016}. We refer to this as the Zaks--Perles--Wills volume conjecture.
\begin{conjecture}\label{ConjectureZPW}
For $d\ge3$ and $k\ge1$, every $d$-dimensional lattice polytope $P$
with exactly $k$ interior lattice points satisfies
\[
\vol(P)\le\vol(S_{d,k})=\frac{(k+1)(s_d-1)^2}{d!}.
\]
With the exception of the case $(d, k) = (3, 1)$, equality holds if and only if $P\simeq S_{d,k}$.
\end{conjecture}

Lagarias and Ziegler \cite{LZ91} gave quantitative volume bounds and proved finiteness, up to affine unimodular equivalence, of lattice polytopes of bounded volume. Consequently, for fixed $d$ and $k \ge 1$, there exist only finitely many affine unimodular equivalence classes. This also lays the groundwork for computational classifications in low dimensions.

\cref{ConjectureZPW} remains open, although some partial results have been obtained.
For simplices with $k=1$, Averkov \cite{Ave12} obtained volume estimates using chains of faces, and Averkov, Kr\"umpelmann, and Nill
\cite{AKN2015} proved the sharp bound.
Balletti, Kasprzyk, and Nill \cite[Corollary~1.2]{BKN22} proved the conjecture for reflexive polytopes in dimension $d\ge4$ and $k=1$, as a consequence of their sharp upper bound for the dual volume of a canonical Fano polytope. 
Averkov, Kr\"umpelmann, and Nill \cite{AKN2020} also obtained a nearly optimal volume bound for all $d$-dimensional lattice simplices with a fixed positive number of interior lattice points.

In the $3$-dimensional case, we have 
\[
S_{3,k}=\conv(0,2e_1,3e_2,6(k+1)e_3),\qquad\vol(S_{3,k})=6(k+1).
\]
Kasprzyk \cite{Kas10} classified $3$-dimensional lattice polytopes with $k=1$; 
in this case there are two extremal affine unimodular equivalence classes. 
Balletti and Kasprzyk \cite[Corollary~5.2]{BK2016} settled the case $k=2$ by classification.
Balletti \cite{Bal21} enumerated all $3$-dimensional lattice simplices with between $1$ and $11$ interior lattice points.
We prove the volume bound and determine the equality case for every $k\ge3$, completing the three-dimensional case. 
The lattice width of $P$, denoted by $\lw(P)$, will be defined in \cref{sec:preliminaries}.

\begin{theorem}\label{main:theorem}
Let $P\subset\R^3$ be a full-dimensional lattice polytope with $k=|\Int(P)\cap\Z^3|\ge3$. Then
\begin{equation}\label{eq:main-bound}
 \vol(P)\le6(k+1).
\end{equation}
Equality holds if and only if $P$ is affinely unimodularly equivalent to $S_{3,k}$. In particular, if $\lw(P)\ge3$, then the inequality is strict.
\end{theorem}

Our proof proceeds by considering the possible lattice widths of $P$.
We first establish the bound $\vol(P)\le6(k+1)$ for polytopes
of width two and show that equality holds precisely for affine
unimodular copies of $S_{3,k}$. We then prove that the inequality
is strict for every larger width. Widths at least $9$ are
handled by a general estimate, leaving only widths $3$ through $8$ to consider. For these widths, the results in the main
text settle all but a few cases involving small numbers of
interior lattice points. These remaining cases are treated
in \cref{s:section}. Together, these arguments establish
both the volume bound and the equality characterization
in \cref{main:theorem}.

The paper is organized as follows.
\Cref{sec:preliminaries} establishes the geometric properties
and area estimates for planar sections used throughout the paper.
\Cref{w2:section} treats width two and determines the equality case.
\Cref{sec:volume} develops the volume estimates based on sections
and projections. In \cref{sec:general}, we combine these results
with those in \cref{s:section} to prove \cref{main:theorem}.
The proof of the sampling inequality is given separately
in \cref{q:proofsection}.

\section{Area bounds for planar sections}\label{sec:preliminaries}

A \emph{convex body} is a compact convex set with nonempty interior; the notation for lattice polytopes extends to convex bodies. For a convex body $K$, set $I(K)=|\Int(K)\cap\Z^d|$ and let $\vol(K)$ be its Euclidean volume. Write $K\simeq L$ if $K$ and $L$ are affinely unimodularly equivalent. Let $\Delta_d=\conv(0,e_1,\dots,e_d)$, where \(e_1,\dots,e_d\) are the standard basis vectors of \(\mathbb R^d\).

For $0\ne u\in\Z^d$, define the \emph{width} in direction $u$ by
\[
\operatorname{width}_u(K)
=\max_{x\in K}\langle u,x\rangle-\min_{x\in K}\langle u,x\rangle.
\]
The \emph{lattice width} is defined as  
\[
\lw(K)=\min_{0\ne u\in\Z^d}\operatorname{width}_u(K).
\]
Here $\langle \cdot, \cdot \rangle$ denotes the standard Euclidean inner product in $\mathbb{R}^d$. The minimum is attained: if $K$ contains a ball of radius $r>0$, then $\operatorname{width}_u(K)\ge2r\|u\|_2$. Thus, only finitely many
integer vectors can have width below any fixed bound. Moreover, a minimizing vector must be primitive. For a full-dimensional lattice polytope, the lattice width is a positive integer. Furthermore, it is at least $2$ whenever the polytope contains an interior lattice point.

For a three-dimensional convex body or polytope $P \subset \mathbb{R}^2 \times \mathbb{R}$, a point in $P$ is conventionally denoted by $(x, z)$, where $x \in \mathbb{R}^2$ and the third coordinate $z \in \mathbb{R}$ represents the \emph{height}.
A \emph{minimum-width slicing} of a lattice three-polytope is an affine unimodular choice of coordinates in which
\begin{equation}\label{eq:normalization}
 P\subset\R^2\times[0,w],\qquad \min_P z=0,\quad\max_P z=w,\quad w=\lw(P).
\end{equation}
For $0\le t\le w$, write
\[ P_t=\{x\in\R^2:(x,t)\in P\},\qquad k_j=I(P_j)\quad(1\le j\le w-1).
\]
Every \emph{section} $P_t$ at height $0<t<w$ is two-dimensional, and its interior is the corresponding section of $\Int(P)$. Hence
\begin{equation}\label{eq:section-accounting}
 k=I(P)=\sum_{j=1}^{w-1}k_j,\qquad \vol(P)=\int_0^w\area(P_t)\,dt,
\end{equation}
where $\area(\cdot)$ denotes Euclidean area.

Since $P$ is convex, for any $a,b\in[0,w]$ and $s\in[0,1]$, we have $(1-s)P_a+sP_b\subseteq P_{(1-s)a+sb}$.
By combining this inclusion with the planar Brunn--Minkowski inequality \cite{Gardner2002}, we obtain that
$t\mapsto\sqrt{\area(P_t)}$ is nonnegative and concave.

We shall also need the mixed-area expansion for Minkowski combinations.
By the two-dimensional case of Minkowski's mixed-volume theorem \cite{Gardner2002}, for nonempty compact convex sets 
$A,B\subset\mathbb{R}^2$ and $s,t\ge 0$,
\begin{equation}\label{eq:mixed-area}
 \area(sA+tB)=s^2\area(A)+2stV(A,B)+t^2\area(B).
\end{equation}
The quantity $V(A,B)$ is called the \emph{mixed area}
(or two-dimensional mixed volume); equivalently,
\[
V(A,B)=\frac{\area(A+B)-\area(A)-\area(B)}{2}.
\]
Mixed area is \emph{additive} under Minkowski addition and \emph{monotone} under inclusion. More precisely, for compact convex sets \(A,B,C\subset\mathbb R^2\), we have 
\[
V(A+B,C)=V(A,C)+V(B,C),
\]
and
\[
A\subseteq B\implies V(A,C)\le V(B,C),
\]
with analogous properties in the second argument. Moreover, for all nonempty compact convex sets \(A,B\subset\mathbb R^2\), we have $V(A,B)^2\ge\operatorname{area}(A)\operatorname{area}(B)$.
If both sets are two-dimensional, equality holds if and only if $A$ and $B$ are \emph{positively homothetic}; that is, $A=x+\lambda B$ for some $x\in\mathbb R^2$ and $\lambda>0$. Degenerate compact
convex sets are included in \eqref{eq:mixed-area} by continuity.
See \cite[Sections~3--5]{Gardner2002} for these classical facts.

\begin{lemma}\label{p:slabs}
Let $P\subset\R^2\times[0,w]$ be a full-dimensional polytope meeting both end planes, where $w\ge2$ is an integer, and suppose
that all its vertices have integer height. Then
\begin{equation}\label{p:slabidentity}
 3\vol(P)=\area(P_0)+\area(P_w)+2\sum_{j=1}^{w-1}\area(P_j)+\sum_{j=0}^{w-1}V(P_j,P_{j+1}).
\end{equation}
Moreover,
\begin{equation}\label{p:mixedpair}
\begin{aligned}
 V(P_{j-1},P_j)+V(P_j,P_{j+1})&\le2\area(P_j) &&(1\le j<w),\\
 V(P_j,P_{j+1})&\ge\sqrt{\area(P_j)\area(P_{j+1})} &&(0\le j<w),
\end{aligned}
\end{equation}
and the end sections satisfy
\begin{small}
\begin{equation}\label{p:endpoints}
\area(P_0)\le\left(2\sqrt{\area(P_1)}-\sqrt{\area(P_2)}\right)^2,\quad 
\area(P_w)\le\left(2\sqrt{\area(P_{w-1})}-\sqrt{\area(P_{w-2})}\right)^2.
\end{equation}
\end{small}
The expressions inside both squares are nonnegative. For brevity, we occasionally write \(V_j:=V(P_j,P_{j+1})\) for \(0\le j<w\), when no confusion arises.
\end{lemma}
\begin{proof}
Between two consecutive integer heights, $P$ has no vertices.
After cutting at these heights, each piece is therefore the convex
hull of its two end sections. Thus
\begin{equation}\label{eq:interval}
P_{j+t}=(1-t)P_j+tP_{j+1}
\qquad(j=0,\ldots,w-1,\ 0\le t\le1).
\end{equation}
Integrating over each unit interval by \eqref{eq:mixed-area}, we have
\[
3\int_j^{j+1}\area(P_t)\,dt
=\area(P_j)+\area(P_{j+1})+V_j.
\]
A summation over all $j$ then proves \eqref{p:slabidentity}.

According to the inclusion $P_{j-1}+P_{j+1}\subseteq2P_j$ together with additivity and monotonicity of mixed area, we have  
\begin{align*}
V(P_{j-1},P_j)+V(P_j,P_{j+1})=V(P_j,P_{j-1}+P_{j+1})\le V(P_j,2P_j)=2\operatorname{area}(P_j).
\end{align*}
This proves the first inequality in \eqref{p:mixedpair}; Minkowski's inequality gives the second.
By the concavity of $t\mapsto\sqrt{\area(P_t)}$, we have \eqref{p:endpoints}, including the asserted nonnegativity. These arguments allow either end section to be a point or a segment. Equation \eqref{eq:interval} shows that $\sqrt{\area(P_t)}$ is continuous with respect to $t$.
\end{proof}

\begin{lemma}\label{p:sections}
Under \eqref{eq:normalization}, for every $0<t<w$, we have 
\[ \lw(P_t)\ge\min\left\{t,w-t\right\}.
\]
\end{lemma}
\begin{proof}
Fix a primitive vector $u\in\Z^2$ and let $Q$ be the image of $P$
under $(x,z)\mapsto(\langle u,x\rangle,z)$.
This integer linear map has rank two, so $Q$ is a two-dimensional
lattice polygon. Every nonzero integer functional on \(Q\) pulls back to a nonzero integer functional on \(P\). Hence $\lw(Q)\ge w$, and its height range $[0,w]$ shows that
$\lw(Q)=w$, attained in the height direction. By \cite[Lemma~2.12]{Bohnert2023}, the horizontal section of $Q$
at height $t$ has length at least
\[
\min\left\{\frac{w}{2},t,w-t\right\}=\min\{t,w-t\}.
\]
This length is exactly $\operatorname{width}_u(P_t)$.
Taking the minimum over primitive integer vectors $u$ proves the claim.
\end{proof}

A convex body is \emph{hollow} if its interior contains no
lattice point. A \emph{split} is a strip
$\{x:m\le\langle u,x\rangle\le m+1\}$, where $u\in\Z^2$ is
primitive and $m\in\Z$. Clearly, any convex body contained in a split is hollow.

\begin{proposition}\label{p:width}
Let $K\subset\R^2$ be a convex body and $n=I(K)$. Set $v=\lw(K)>1$. If $1<v\le2$, then
\[\area(K)\le\frac{v^2}{2(v-1)}(n+1).
\]
If $v\ge2$, then
\[\area(K)\le\begin{cases}2,&n=0,\\ 2n+\frac{5}{2},&n\ge1.\end{cases}
\]
In particular, in this range $\area(K)\le\frac94(n+1)$.
\end{proposition}
\begin{proof}
First assume that $n\ge1$.  The first bound and the bound for $2<v\le5$ follow from \cite[Theorem~1.2]{Bohnert2023}; 
for $v>5$ use the stronger bound in \cite[Theorem~1.1]{Bohnert2023}. 
At $v=2$, the first bound equals $2(n+1)\le2n+\frac52$.
Hence the assertion for $n\ge1$ follows.

Now assume that $n=0$. For $1<v\le2$, the first bound is proved in
\cite[Theorem~2.2]{AW2012}; at $v=2$ it gives
\[
\area(K)\le\frac{v^2}{2(v-1)}(n+1)=2.
\]
For $v>2$, the results in \cite[Theorems~2.1 and 2.2]{AW2012} show $v\le1+\frac2{\sqrt3}$ and
\[
\area(K)\le\frac{3v^2}{3v+1-\sqrt{1+6v-3v^2}}\le2.
\]
Indeed,
\[
2\left(3v+1-\sqrt{1+6v-3v^2}\right)-3v^2
=\left(\sqrt{1+6v-3v^2}-1\right)^2\ge0.
\]
Finally, $2n+\frac52\le\frac94(n+1)$ for $n\ge1$, and $2\le\frac94$
for $n=0$; hence $\area(K)\le\frac94(n+1)$ whenever $v\ge2$.
\end{proof}

For a positive integer $D$ define
\[
 \mathcal G_D=\bigcup_{d=1}^D d^{-1}\Z^2, \qquad A_D=\frac{(D+1)^2}{2D}.
\]
Integer-height sections of lattice polytopes have vertices in $\mathcal G_D$ for an appropriate $D$. 
The vertices may have different denominators, whose least common multiple can exceed $D$.

\begin{lemma}\label{p:denominator}
Let $K$ be a two-dimensional convex polygon with vertices in $\mathcal G_D$, where $D\ge2$. Suppose that $K$ is not contained
in a split. Set $n=I(K)$, and assume that either $D\ge3$ or $n\ge2$. Then
\[
 \area(K)\le A_D(n+1).
\]
Equality holds if and only if $K$ is affinely unimodularly equivalent to
\[
 E_{D,n}=\conv\left\{(0,0),\left(1+\frac{1}{D},0\right),(0,(D+1)(n+1))\right\}.
\]
In particular, equality implies $\lw(K)=1+\frac{1}{D}$.
\end{lemma}
\begin{proof} 
After a unimodular change of coordinates, we may assume that the first coordinate realizes the minimum lattice width of \(K\).  An integer translation makes the first integer in the interior of its range equal to zero. Such an integer exists, since otherwise $K$ would lie in a split. Thus the range has the form $[-a,N-1+b]$, where $N\ge1$ is the number of interior integer coordinates and  
\[D^{-1}\le a,b\le1.
\]
Here $0<a,b\le1$ because the adjacent integers lie outside the interior of the projection interval. Its endpoints are coordinates
of vertices in $\mathcal G_D$. Therefore, we have $a,b\ge D^{-1}$.

Suppose first that $N=1$. By reflecting the first coordinate if necessary, we may assume that $a\le b$. Let $h(x)$ be the length of the vertical section of $K$. The convexity of \(K\) implies that \(h\) is concave. All interior lattice points lie on $x=0$. Thus $h(0)\le n+1$. Since $h$ is concave, there exists \(s\in\mathbb R\) such that $h(x)\le h(0)+sx$ for all \(x\) in the domain of \(h\). Its nonnegativity at $-a$ implies $s\le \frac{h(0)}{a}$.
Integrating the inequality $h(x)\le h(0)+sx$ over $[-a,b]$ and using $b\ge a$, we obtain
\[
 \area(K)\le h(0)(a+b)+\frac12s(b^2-a^2)\le\frac{(a+b)^2}{2a}h(0)
 \le\frac{(1+a)^2}{2a}(n+1)\le A_D(n+1).
\]
The last inequality follows because $\frac{(1+a)^2}{2a}$ decreases on $(0,1]$. 
Equality requires $a=\frac{1}{D}$ and $b=1$, hence $\lw(K)=1+\frac{1}{D}$.

If $N\ge2$ and $\lw(K)\le2$, then $\lw(K)=N-1+a+b\ge1+\frac{2}{D}$. The function $\frac{v^2}{2(v-1)}$ decreases
on $(1,2]$. Therefore, by \cref{p:width}, we have
\[
\area(K)\le\frac{(D+2)^2}{4D}(n+1)<A_D(n+1).
\]
Suppose that $\operatorname{lw}(K)>2$.
If $D\ge 3$, by applying \cref{p:width} again, we obtain 
\[
\operatorname{area}(K)\le \frac{9}{4}(n+1) < A_D(n+1).
\]
If $D=2$, then $n\ge 2$ by assumption, and we obtain
\[
\operatorname{area}(K)
\le 2n+\frac{5}{2}
< \frac{9}{4}(n+1)
= A_2(n+1).
\]

It remains to identify equality in the first case. Equality forces
$a=\frac{1}{D}$, $b=1$, $h(0)=n+1$, and $h(x)=(n+1)(1+Dx)$.
The lower boundary of $K$ is convex and the upper boundary is concave; since their difference is affine, both boundaries are affine. Hence \(K\) is a triangle with one vertex at \(\left(-\frac1D,\alpha\right)\) and the opposite edge lying on the line \(x=1\), of length \((D+1)(n+1)\). The vertex $\left(-\frac{1}{D},\alpha\right)$ belongs to $\mathcal G_D$. Therefore, its first coordinate forces $D\alpha\in\Z$.

The open section at $x=0$ has length $n+1$ and contains exactly $n$ integers. It follows that its lower endpoint is an integer $q$. The lower endpoint of the vertical edge at \(x=1\) therefore has second coordinate $(D+1)q-D\alpha$. The affine unimodular map
\[
 (x,y)\longmapsto\bigl(x,y+(D\alpha-Dq)x-q\bigr)
\]
sends the three vertices to $(-\frac{1}{D},0)$, $(1,0)$, and $(1,(D+1)(n+1))$. The reflection $x\mapsto1-x$ maps the polytope onto $E_{D,n}$.

Conversely, $E_{D,n}$ has exactly $n$ interior lattice points, area $A_D(n+1)$, and lattice width $1+1/D>1$. Thus $E_{D,n}$ is not contained in a split. Since these properties are preserved under affine unimodular transformations, the same holds for every polygon affinely unimodularly equivalent to $E_{D,n}$.
\end{proof}

The rectangles $[0,1]\times[0,N]$ ($N\in\Z_{>0}$) are hollow and have arbitrarily large area. Thus the split exclusion in \cref{p:denominator} is necessary. 

\begin{lemma}\label{p:rational}
Let $P$ be a lattice three-polytope with height range $[0,w]$. Then every vertex of an internal integer-height section
belongs to $\mathcal G_w$. Consequently, for $w\ge3$,
a section $P_j$ not contained in a split satisfies
\[
\area(P_j)\le A_w(k_j+1).
\]
Equality in this bound implies
$\lw(P_j)=1+\frac1w$.
\end{lemma}
\begin{proof}
A section vertex is either a vertex of $P$ or the intersection of
an edge with the section plane. For integer edge endpoints
$(x,r),(y,s)$ with $r<s$ and $r<j<s$, the intersection has
horizontal coordinate
\[
\frac{(s-j)x+(j-r)y}{s-r}\in\frac1{s-r}\Z^2,
\qquad s-r\le w.
\]
This proves the denominator assertion.
The area bound and its equality condition follow
from \cref{p:denominator}.
\end{proof}

Finally, we mention that estimates for the area of rational polygons have also been studied in~\cite{BS2024}.

\section{Width two and the equality case}\label{w2:section}

At width two, every interior lattice point lies in the middle section. 
The volume bound follows by combining a sharp area estimate for this half-integral section with a sharp comparison between its
area and the volume of the polytope. 
The equality conditions identify the section and force the polytope to be a pyramid.

\begin{lemma}\label{w2:halfarea}
Let $Q\subset\R^2$ be a two-dimensional convex polygon such that $2Q$ is a lattice polygon, and suppose that $i=|\Int(Q)\cap\Z^2|\ge2$.  Then
\begin{equation}\label{w2:area}
 \area(Q)\le\frac94(i+1).
\end{equation}
Equality holds precisely when $Q$ is affinely unimodularly equivalent to
\[
 T_i=\conv\left\{\left(0,\frac12\right),(0,-1),(3(i+1),-1)\right\}.
\]
\end{lemma}
\begin{proof}
The vertices of $Q$ belong to $\mathcal G_2$, and its interior lattice
points prevent it from being contained in a split. We consider \cref{p:denominator} with $D=2$. Its equality triangle $E_{2,i}$
is carried to $T_i$ by the affine unimodular map $(x,y)\mapsto(y,x-1)$.
\end{proof}

\begin{lemma}\label{w2:middle}
Let $P\subset\R^2\times[0,2]$ be a three-dimensional convex polytope meeting both end planes.
Then
\begin{equation}\label{w2:middlebound}
 \vol(P)\le\frac83\area(P_1).
\end{equation}
Equality holds if and only if, after possibly interchanging the two end planes, $P_0=\{v\}$ and $P_2=2P_1-v$. 
In this event
\[
 P=\conv\bigl(\{(v,0)\}\cup((2P_1-v)\times\{2\})\bigr).
\]
\end{lemma}
\begin{proof}
The function $f(t)=\sqrt{\area(P_t)}$ is continuous, nonnegative,
and concave. Thus there exists $s\in\R$ such that
$f(t)\le f(1)+s(t-1)$ for $0\le t\le2$.
The right-hand side is nonnegative at both endpoints, so $|s|\le f(1)$.
Consequently,
\begin{align*}
 \vol(P)=\int_0^2 f(t)^2\,dt\le\int_0^2\left(f(1)+s(t-1)\right)^2\,dt
 =2f(1)^2+\frac23s^2\le\frac83f(1)^2=\frac83\area(P_1).
\end{align*}
Equality forces $s=\pm f(1)$ and $f(t)=f(1)+s(t-1)$.
After reflecting height if necessary, $f(t)=ct$ for some $c>0$.
Let $v\in P_0$ and set
\[
 C=\conv\bigl(\{(v,0)\}\cup(P_2\times\{2\})\bigr)\subseteq P.
\]
For every $0<t<2$, the two-dimensional sections satisfy
$C_t\subseteq P_t$ and
\[
 \area(C_t)=\left(\frac{t}{2}\right)^2\area(P_2)=c^2t^2=\area(P_t).
\]
Thus $C_t=P_t$. Furthermore, by taking closures, we get $C=P$.
It follows that $P_0=\{v\}$ and $P_1=\frac{v+P_2}{2}$, as required.
Conversely, by the stated endpoint conditions, we have $f(0)=0$ and $f(2)=2f(1)$. 
The concavity of $f(t)$ then forces $f(t)=t f(1)$. The same section comparison yields the stated pyramid, whose volume is \(\frac83\operatorname{area}(P_1)\).

This completes the proof.
\end{proof}

\begin{remark}\label{rem:midareatovol}
Let $P$ be a three-dimensional convex polytope with height range $[a,b]$, where $a<b$. We first apply the affine transformation $z\mapsto\frac{2(z-a)}{b-a}$ to the height coordinate, and then invoke \cref{w2:middle} to obtain
\begin{equation}\label{w2:rescaled-middle}
 \vol(P)\le\frac{4(b-a)}{3}\area\left(P_{\frac{a+b}{2}}\right).
\end{equation}
Equality holds precisely for a pyramid with its apex in one end plane
and its base in the other.
\end{remark}

\begin{theorem}\label{w2:theorem}
Let $P\subset\R^3$ be a lattice polytope with $\lw(P)=2$ and $k=|\Int(P)\cap\Z^3|\ge2$.  Then
\[
 \vol(P)\le6(k+1),
\]
with equality if and only if $P$ is affinely unimodularly equivalent to
\[
 S_{3,k}=\conv\{0,2e_1,3e_2,6(k+1)e_3\}.
\]
\end{theorem}
\begin{proof}
We consider the coordinates as in \eqref{eq:normalization} with $w=2$.
A vertex of the section $P_1$ is either an integral vertex of $P$ at that height or the midpoint of an edge joining the two end planes.  Thus $2P_1$ is integral. By \eqref{eq:section-accounting}, we have $I(P_1)=k$.
By \cref{w2:middle,w2:halfarea}, we obtain 
\[
 \vol(P)\le\frac83\area(P_1) \le\frac83\cdot \frac94(k+1)=6(k+1).
\]

Suppose equality holds.  By \cref{w2:middle}, after a reflection in height if necessary, $P$ is a pyramid with integral apex $(v,0)$
and base $(2P_1-v)\times\{2\}$.  The affine unimodular shear
\[
 (x,z)\longmapsto(x+(z-1)v,z)
\]
fixes the plane $z=1$ pointwise, sends the apex to the origin, and sends the base to $2P_1\times\{2\}$.  A planar affine unimodular map $x\mapsto Ux+t$, with $U\in\mathrm{GL}_2(\Z)$ and $t\in\Z^2$, lifts to the three-dimensional unimodular map
\[
 (x,z)\longmapsto(Ux+zt,z).
\]
It fixes the apex and induces the required map on the middle section. \cref{w2:halfarea} therefore allows us to assume $P_1=T_k$.

The resulting pyramid has vertices
\[
 (0,0,0),\quad(0,1,2),\quad(0,-2,2),\quad(6(k+1),-2,2).
\]
After translating $(0,-2,2)$ to the origin, its other vertices are
\[
 2(0,1,-1),\qquad3(0,1,0),\qquad6(k+1)(1,0,0).
\]
The three primitive vectors displayed here form a unimodular basis.
Hence this tetrahedron is equivalent to $S_{3,k}$.

Conversely, the interior lattice points of $S_{3,k}$ are precisely $(1,1,j)$ for $1\le j\le k$: its strict facet inequalities force the first two coordinates to be one and the last to lie in this range.
Its volume is $6(k+1)$.  Its first coordinate gives width at most two, while an interior lattice point rules out width one.  Thus it
has width two and realizes equality.
\end{proof}

\section{Volume bounds from sections and projections}\label{sec:volume}
We derive volume bounds in two ways: from the areas of integer-height
sections, and from the lengths of fibers over a planar projection.
These will be applied according to whether an internal integer-height
section is contained in a split.

\subsection{A sharp sampling inequality}\label{q:section}
To identify the optimal sampling constant, consider a pyramid with apex at height zero and base at an integer height $w\ge2$. Its section at height $t$ is a translate of the base scaled by $\frac{t}{w}$. Thus its area is proportional to \(t^2\). This motivates the definition
\[
C_w := \frac{\int_0^w t^2 \, dt}{\sum_{j=1}^{w-1} j^2} = \frac{2w^2}{(w-1)(2w-1)}.
\]
Any constant in such a bound must be at least $C_w$. The following theorem shows that $C_w$ also suffices for products of two nonnegative concave functions. The proof in \cref{q:proofsection} is based on \cite[Theorem~4.1]{HM1995}.

\begin{theorem}\label{q:bilinear}
Let $w\ge2$ be an integer. For finite nonnegative concave functions $f,g:[0,w]\to\R$,
\[
 \int_0^w f(t)g(t)\,dt \le C_w\sum_{j=1}^{w-1}f(j)g(j).
\]
The constant $C_w$ is optimal.
\end{theorem}

By taking $f(t)=g(t)=\sqrt{\area(P_t)}$ in \cref{q:bilinear} and using \eqref{eq:section-accounting},
we obtain the following corollary. 
\begin{corollary}\label{q:volume}
Let $P\subset\R^2\times[0,w]$ be a three-dimensional compact convex body meeting both end planes, where $w\ge2$ is an integer.
For $P_t=\{x\in\R^2:(x,t)\in P\}$, we have 
\begin{equation}\label{q:volume-ineq}
 \vol(P)\le C_w\sum_{j=1}^{w-1}\area(P_j).
\end{equation}
The constant $C_w$ is optimal, as shown by pyramids with an apex on one end plane and a two-dimensional base on the other.
\end{corollary}

\subsection{Split sections and triangular projections}\label{sec:projection}

In a minimum-width slicing, a section contained in a split determines the projection described below. 
Its planar classification is a special case of \cite[Lemma~4.3]{Bohnert2023}. 
A \emph{surjective lattice projection} is an affine map carrying $\Z^3$ onto $\Z^2$.

\begin{lemma}\label{pr:projection}
Suppose that \eqref{eq:normalization} holds with $w\ge3$. If an internal integer-height section $P_j$ is contained in a split, then
$j\in\{1,w-1\}$ and $P$ admits a surjective lattice projection onto a polygon affinely unimodularly equivalent to $w\Delta_2$.
\end{lemma}

\begin{proof}
Write the split as $m\le\langle u,x\rangle\le m+1$,
where $u\in\Z^2$ is primitive and $m\in\Z$. Then
$\pi(x,z)=(\langle u,x\rangle-m,z)$ is a surjective
lattice projection. As in the proof of \cref{p:sections},
$Q=\pi(P)$ is a lattice polygon with minimum lattice
width $w$ and height range $[0,w]$.

By \cref{p:sections}, we have $\min\{j,w-j\}\le\lw(P_j)\le1$, so $j\in\{1,w-1\}$.
Moreover, $Q\cap(\R\times\{j\})\subseteq[0,1]\times\{j\}$.
Thus the line $z=j$ meets $\Int(Q)$ but contains no interior lattice point of $Q$. 
After interchanging the coordinates, we apply the result in \cite[Lemma~4.3]{Bohnert2023}. Since the height direction realizes $\lw(Q)$, this shows that $Q$ is affinely unimodularly equivalent to an integer multiple of $\Delta_2$.
That multiple is $w$, since $\lw(Q)=w$.
This completes the proof.
\end{proof}

To estimate volume from the projection, we use the following sampling lemma. Here, recall that a finite concave function on a triangle is bounded and continuous in the interior, hence integrable.

\begin{lemma}\label{pr:sampling}
Let $w\ge3$ be an integer. Every finite concave function
$h:w\Delta_2\to\R$ satisfies
\begin{equation}\label{eq:triangle-sampling}
 \int_{w\Delta_2}h(x)\,dx\le \frac{w^2}{(w-1)(w-2)} \sum_{p\in\Int(w\Delta_2)\cap\Z^2}h(p).
\end{equation}
For nonnegative functions, the multiplicative constant is optimal.
\end{lemma}

\begin{proof}
Set $v_0=(0,0)$, $v_1=(w,0)$, and $v_2=(0,w)$.
Let $\lambda_i(x)$ be the barycentric coordinates. Thus we have 
\[
 x=\sum_{i=0}^2\lambda_i(x)v_i, \qquad \lambda_i(x)\ge0, \qquad \sum_{i=0}^2\lambda_i(x)=1.
\]
Set $\mathcal A=\left\{\alpha\in\Z_{\ge0}^3:|\alpha|=w-3\right\}$,
where $|\alpha|=\alpha_0+\alpha_1+\alpha_2$. Define
\[
 B_\alpha(x)=
 \frac{(w-3)!}{\alpha_0!\alpha_1!\alpha_2!}
 \lambda_0(x)^{\alpha_0}
 \lambda_1(x)^{\alpha_1}
 \lambda_2(x)^{\alpha_2}.
\]
These weights are nonnegative, and the multinomial theorem shows $\sum_{\alpha\in\mathcal A}B_\alpha=1$.

For nonnegative integers $\beta_i$, set $|\beta|=\beta_0+\beta_1+\beta_2$.
By repeated integration, we get 
\[
\int_{w\Delta_2}\lambda_0^{\beta_0}\lambda_1^{\beta_1}
\lambda_2^{\beta_2}\,dx
=w^2\frac{\beta_0!\beta_1!\beta_2!}{(|\beta|+2)!}.
\]
By applying this formula to $B_\alpha$ and $\lambda_iB_\alpha$, we obtain 
\[
 \int_{w\Delta_2} B_\alpha(x)\,dx
 =\frac{w^2}{(w-1)(w-2)}=:M,
 \qquad
 \frac1M\int_{w\Delta_2}\lambda_i(x)B_\alpha(x)\,dx
 =\frac{\alpha_i+1}{w}
 \quad(i=0,1,2).
\]
Since $x_1=w\lambda_1(x)$ and $x_2=w\lambda_2(x)$,
the density $B_\alpha$ has barycenter
\[
 p_\alpha
 :=\frac1M\int_{w\Delta_2} xB_\alpha(x)\,dx
 =(\alpha_1+1,\alpha_2+1).
\]

Both coordinates of $p_\alpha$ are positive integers, and
\[
 (\alpha_1+1)+(\alpha_2+1)=w-1-\alpha_0<w.
\]
Hence $p_\alpha\in\Int(w\Delta_2)\cap\Z^2$.
Conversely, every interior lattice point $(r,s)$ corresponds
to the unique index
$\alpha=(w-1-r-s,r-1,s-1)\in\mathcal A$.
Thus the points $p_\alpha$ exhaust $\Int(w\Delta_2)\cap\Z^2$
without repetition.

By applying Jensen's inequality to the normalized density $B_\alpha/M$, we obtain
\[
 \frac1M\int_{w\Delta_2} B_\alpha(x)h(x)\,dx
 \le
 h\left(\frac1M\int_{w\Delta_2} xB_\alpha(x)\,dx\right)
 =h(p_\alpha).
\]
By multiplying the preceding inequality by $M$ and summing over $\mathcal A$, we obtain \eqref{eq:triangle-sampling}.

To prove optimality, take $h\equiv1$, which is nonnegative
and concave. A direct calculation gives
\[
 \area(w\Delta_2)=\frac{w^2}{2},
 \qquad
 I(w\Delta_2)=\frac{(w-1)(w-2)}2.
\]
Hence any admissible constant $C$ must satisfy
\[
 C\ge\frac{\area(w\Delta_2)}{I(w\Delta_2)} =\frac{w^2}{(w-1)(w-2)}.
\]
This proves optimality. 
\end{proof}

\begin{proposition}\label{pr:volume}
Let $P$ be a full-dimensional lattice three-polytope with $I(P)=k$.
If a surjective lattice projection maps $P$ onto a polygon equivalent to $w\Delta_2$, where $w\ge3$ is an integer, then
\begin{equation}\label{eq:projection-volume}
 \vol(P)\le\frac{w^2}{(w-1)(w-2)}k+\frac{w^2}{2}.
\end{equation}
\end{proposition}
\begin{proof}
Let $A$ denote the linear part of the projection. The surjectivity of $A:\mathbb Z^3\to\mathbb Z^2$ ensures that any lifts of a lattice basis of $\mathbb Z^2$ can be extended to a lattice basis of $\mathbb Z^3$ by adjoining a primitive generator of $\ker A\cap\mathbb Z^3$.
Thus, after affine unimodular changes of coordinates,
we may assume that
\[
 \pi(x,z)=x, \qquad \pi(P)=w\Delta_2.
\]
For $x\in w\Delta_2$, write
\[
 P\cap(\{x\}\times\R)
 =\{x\}\times[\ell(x),u(x)],
 \qquad h(x):=u(x)-\ell(x).
\]
Obviously $\vol(P)=\int_{w\Delta_2} h(x)\,dx$. By the convexity of $P$, the function $u$ is concave and $\ell$ is convex. Since $u\ge\ell$, their difference $h=u-\ell$ is nonnegative and concave.

For every $p\in\Int(w\Delta_2)\cap\Z^2$, we have
\[
 \Int(P)\cap(\{p\}\times\R) =\{p\}\times(\ell(p),u(p)).
\]
Define
$n(p):=\bigl|(\ell(p),u(p))\cap\Z\bigr|.$
An open interval containing $n(p)$ integers has length
at most $n(p)+1$, hence
$h(p)\le n(p)+1$.
Every interior lattice point of $P$ projects to a point of $\Int(w\Delta_2)\cap\Z^2$. Thus 
\[
 \sum_{p\in\Int(w\Delta_2)\cap\Z^2}n(p)=k,
 \qquad
 |\Int(w\Delta_2)\cap\Z^2|=I(w\Delta_2)=\frac{(w-1)(w-2)}2.
\]
Consequently,
\[
 \sum_{p\in\Int(w\Delta_2)\cap\Z^2}h(p)
 \le \sum_{p\in\Int(w\Delta_2)\cap\Z^2}\bigl(n(p)+1\bigr)
 =k+\frac{(w-1)(w-2)}2.
\]
The conclusion follows from \cref{pr:sampling}.
\end{proof}

\begin{remark}
The bound \eqref{eq:projection-volume} is sharp: equality holds for $P=w\Delta_2\times[0,L]$,
where $L\ge1$ is an integer.
Moreover, $\lw(P)=\min\{w,L\}$.
\end{remark}

By combining \cref{pr:projection,pr:volume}, we see that
\eqref{eq:projection-volume} applies whenever a minimum-width
slicing with $w=\lw(P)\ge3$ has an internal integer-height
section contained in a split.

\section{Strict volume bounds for larger lattice widths}\label{sec:general}

We prove that \(\operatorname{vol}(P)<6(k+1)\) whenever \(k\ge3\) and \(w=\operatorname{lw}(P)\ge3\). We first establish this inequality for \(w\ge9\), leaving only widths $3$ through $8$. For \(w=3\), we combine the area and volume bounds and show that equality cannot occur. For \(4\le w\le8\), the estimates below reduce the proof to the cases \(4\le w\le7\) and \(3\le k<w\), which are treated in \cref{s:section}.

\subsection{Large widths and small interior-point counts}

Let $\operatorname{Flt}(d)$ denote the supremum of the lattice widths of hollow convex bodies in $\R^d$. 
The estimate of Averkov, Codenotti, Macchia and Santos \cite[Theorem~5.2]{ACMS2021} implies
\begin{equation}\label{l:flatness}
 \operatorname{Flt}(3)<3.972<4.
\end{equation}

We combine the estimate of Freyer and Lucas \cite[Proposition~1.6]{FL2022} with a covering condition obtained from the definition of \(\operatorname{Flt}(d)\).

\begin{lemma}\label{l:count}
  Let \(K\) be a convex body in \(\mathbb R^d\), and let \(0<t<1\). If
\[
t\operatorname{lw}(K)>\operatorname{Flt}(d),
\]
then \(tK+\mathbb Z^d=\mathbb R^d\) and
\[
I(K)\ge(1-t)^d\operatorname{vol}(K).
\]
\end{lemma} 

\begin{proof} 
Suppose that \(tK+\mathbb Z^d\ne\mathbb R^d\). Let \(y\notin tK+\mathbb Z^d\). Then \(y-tK\) contains no lattice point. By the definition of \(\operatorname{Flt}(d)\), we have 
\[
t\operatorname{lw}(K)=\operatorname{lw}(y-tK)\le\operatorname{Flt}(d),
\]
a contradiction. Hence \(tK+\mathbb Z^d=\mathbb R^d\). Furthermore, the covering radius $\mu(K):=\inf\{s>0:sK+\Z^d=\R^d\}$ satisfies $\mu(K)\le t<1$. By \cite[Proposition~1.6]{FL2022}, we obtain 
$$I(K)\ge(1-\mu(K))^d\vol(K)\ge(1-t)^d\vol(K).$$
This completes the proof.
\end{proof}

\begin{proposition}\label{l:large}
If a lattice three-polytope $P$ has width at least $9$ and $k\ge1$ interior lattice points, then
\[
 \vol(P)\le\frac{729}{125}k<6(k+1).
\]
\end{proposition}
\begin{proof}
Let $t=\frac{4}{9}$. Then $t\lw(P)\ge4>\operatorname{Flt}(3)$. By \cref{l:count}, we have $k\ge\left(\frac{5}{9}\right)^3\vol(P)$. Since $\frac{729}{125}<6$, we complete the proof.
\end{proof}

\begin{proposition}\label{s:parity}
A lattice three-polytope with at most $7$ interior lattice points has lattice width at most $7$.
\end{proposition}
\begin{proof}
Let \(P\) be such a polytope. One of the eight residue classes modulo $2\Z^3$, say $a+2\Z^3$, contains no interior lattice point. Thus $\frac{P-a}{2}$ is hollow,
and the flatness estimate \eqref{l:flatness} gives $\frac{\lw(P)}{2}<4$.
The lattice width of a lattice polytope is an integer, since its width in every integer direction is an integer and the minimum is
attained. Consequently $\lw(P)\le7$.
\end{proof}

\subsection{Width three}\label{w3:section}

We now consider the case \(w=3\). Since all interior lattice points of \(P\) lie in the two sections \(P_1\) and \(P_2\), we use these sections to prove that \(\operatorname{vol}(P)<6(k+1)\) for every \(k\ge3\).

\begin{lemma}\label{w3:geometry}
Let $P\subset\R^2\times[0,3]$ be a three-dimensional polytope meeting both end planes, with every vertex at integer height.
Let $q=\area(P_1)$ and $r=\area(P_2)$. Then
\[
 \vol(P)\le3\bigl(q+r-\sqrt{qr}\bigr)\le\frac95(q+r).
\]
If $\vol(P)=\frac{9}{5}(q+r)$, then $P_1$ and $P_2$ are positively homothetic and
\begin{equation}\label{w3:equality-ratio}
 \frac{q}{r}\in\left\{\frac{1}{4},4\right\}.
\end{equation}
\end{lemma}
\begin{proof}
By \cref{p:slabs}, we have 
\begin{align*}
 V_{j-1}+V_{j}\le2\area(P_j) (j=1,2),\qquad  V_1\ge\sqrt{qr},
\end{align*}
and 
\begin{equation*}
\area(P_0)\le\bigl(2\sqrt{q}-\sqrt{r}\bigr)^2,\quad 
\area(P_3)\le\bigl(2\sqrt{r}-\sqrt{q}\bigr)^2.
\end{equation*}
By substituting these inequalities into \eqref{p:slabidentity}, we obtain
\begin{align*}
  3\vol(P)=&\area(P_0)+\area(P_3)+2q+2r+V_0+V_1+V_2\\
  \le&\bigl(2\sqrt{q}-\sqrt{r}\bigr)^2+\bigl(2\sqrt{r}-\sqrt{q}\bigr)^2+2q+2r+(V_0+V_1)+(V_1+V_2)-V_1\\
  \le&\bigl(2\sqrt{q}-\sqrt{r}\bigr)^2+\bigl(2\sqrt{r}-\sqrt{q}\bigr)^2+4q+4r-\sqrt{qr}\\
  =&9\bigl(q+r-\sqrt{qr}\bigr).
\end{align*}
According to the nonnegative concavity of $\sqrt{\area(P_t)}$, we have $\frac{\sqrt{q}}{2}\le \sqrt{r}\le2\sqrt{q}$ and 
\[
 5\sqrt{qr}-2(q+r)=(2\sqrt{q}-\sqrt{r})(2\sqrt{r}-\sqrt{q})\ge0.
\]
This proves both bounds. Equality in the final bound forces
$V_1=\sqrt{qr}$ and $\frac{q}{r}\in\left\{\frac{1}{4},4\right\}$. Equality in Minkowski's inequality
makes $P_1,P_2$ positively homothetic, and their area ratio is
$\frac{1}{4}$ or $4$.
\end{proof}

\begin{theorem}\label{w3:strict}
Every $3$-dimensional lattice polytope $P$ of lattice width $3$ with $k=|\Int(P)\cap\Z^3|\ge3$ satisfies
\[
 \vol(P)<6(k+1).
\]
\end{theorem}
\begin{proof}
We consider the coordinates as in \eqref{eq:normalization} with $w=3$. Let $q=\area(P_1)$ and $r=\area(P_2)$.
By \eqref{eq:section-accounting}, we have $k_1+k_2=k$.
If either $P_1$ or $P_2$ is contained in a lattice split, \cref{pr:projection} gives a surjective lattice projection of $P$ onto a polygon unimodularly equivalent to $3\Delta_2$. By \cref{pr:volume} with $w=3$, it follows that 
\[
 \vol(P)\le\frac92(k+1)<6(k+1).
\]

Suppose now that neither internal section is contained in a split.
We apply \cref{p:rational} to both sections, including any hollow section, to obtain
\[
 q\le\frac83(k_1+1),\qquad
 r\le\frac83(k_2+1).
\]
\cref{w3:geometry} now implies
\begin{equation}\label{w3:count-bound}
 \vol(P)\le\frac95(q+r)\le\frac{24}{5}(k+2).
\end{equation}
For $k>3$, the last expression is strictly less than $6(k+1)$.

It remains to exclude equality for $k=3$. If $\vol(P)=24$, then every inequality in \eqref{w3:count-bound} is an equality. 
We recall the equality condition in \cref{p:rational}. It follows that 
\[
 \lw(P_1)=\lw(P_2)=\frac43.
\]
On the other hand, \cref{w3:geometry} makes $P_1,P_2$ positively homothetic with area ratio $4$ or $\frac{1}{4}$. Their
homothety ratio, and hence their lattice-width ratio, is therefore $2$ or $\frac{1}{2}$. This contradicts their equal lattice widths. Thus equality is impossible.
\end{proof}

\subsection{Widths four through eight}

For $w\ge4$, the sections $P_2,\ldots,P_{w-2}$ have width at least two by \cref{p:sections}. 
Therefore, we may consider using \cref{p:width}.
The following estimate bounds their share of the total sampled area.

\begin{lemma}\label{g:share}
For an integer $w\ge4$ and a nonnegative concave function $f:[0,w]\to\R$, we have 
\[
\sum_{j=2}^{w-2}f(j)^2\ge r_w\sum_{j=1}^{w-1}f(j)^2,\qquad
r_w=1-\frac{6\bigl(1+(w-1)^2\bigr)}{w(w-1)(2w-1)}.
\]
The constant is optimal. If $f$ is nonzero and equality holds, then $\frac{f(1)}{f(w-1)}\in\left\{w-1,\frac{1}{w-1}\right\}$.
\end{lemma}
\begin{proof}
If $f=0$, the result is immediate. Otherwise, nonnegative concavity of $f$ gives that $a=f(1)>0$, $b=f(w-1)>0$ and
$a\ge \frac{b}{w-1}$, $b\ge \frac{a}{w-1}$. By expanding $(a-(w-1)b)(b-(w-1)a)\ge0$, we get 
\begin{equation}\label{eq:end-sample-product}
 \frac{ab}{a^2+b^2}\ge\frac{w-1}{1+(w-1)^2},
\end{equation}
where equality holds if and only if $\frac{a}{b}\in\left\{w-1,\frac{1}{w-1}\right\}$. By the concavity of $f$, we have
\[
\begin{aligned}
\sum_{j=1}^{w-1}f(j)^2
&\ge\sum_{j=1}^{w-1}
\left(\frac{(w-1-j)a+(j-1)b}{w-2}\right)^2\\
&=\frac{w-1}{6(w-2)}
\bigl((2w-3)(a^2+b^2)+2(w-3)ab\bigr)\\
&\ge\frac{w(w-1)(2w-1)}
{6\bigl(1+(w-1)^2\bigr)}(a^2+b^2),
\end{aligned}
\]
where the last inequality follows from \eqref{eq:end-sample-product}. Then, we obtain
\[
\begin{aligned}
\sum_{j=2}^{w-2}f(j)^2=\sum_{j=1}^{w-1}f(j)^2-(a^2+b^2)\ge r_w\sum_{j=1}^{w-1}f(j)^2.
\end{aligned}
\]
Equality can hold only if equality holds in \eqref{eq:end-sample-product}.
The functions $f(t)=t$ and $f(t)=w-t$ show optimality.
\end{proof}

\begin{proposition}\label{g:nonsplit}
Suppose that \(P\) satisfies \eqref{eq:normalization} with \(w\ge4\), and neither \(P_1\) nor \(P_{w-1}\) is contained in a split. Then
\[
\operatorname{vol}(P)\le B_w(k+w-1),
\qquad
B_w=\frac{C_w}{\frac{1-r_w}{A_w}+\frac{4r_w}{9}}.
\]
If equality holds, then
\[
\frac{\operatorname{area}(P_1)}{\operatorname{area}(P_{w-1})}\in\left\{(w-1)^2,\frac{1}{(w-1)^2}\right\},
\]
and
\[
\operatorname{area}(P_j)=A_w(k_j+1)
\qquad(j=1,w-1).
\]
\end{proposition}
\begin{proof}
By \cref{p:sections}, the sections \(P_j\) with \(2\le j\le w-2\) have lattice width at least two. 
By \cref{p:width,p:rational}, we obtain
\[
\operatorname{area}(P_j)\le
\begin{cases}
A_w(k_j+1),&j=1,w-1,\\[2pt]
\dfrac94(k_j+1),&2\le j\le w-2.
\end{cases}
\]
Set $S=\sum_{j=1}^{w-1}\operatorname{area}(P_j)$. By applying \cref{g:share} to the nonnegative concave function
\(f(t)=\sqrt{\operatorname{area}(P_t)}\), we obtain
\[
\sum_{j=2}^{w-2}\operatorname{area}(P_j)\ge r_wS.
\]
Since \(A_w>\frac{9}{4}\), it follows that
\begin{align*}
k+w-1=\sum_{j=1}^{w-1}(k_j+1)\ge\frac{S}{A_w}+\left(\frac49-\frac1{A_w}\right)
 \sum_{j=2}^{w-2}\operatorname{area}(P_j)\ge\left(\frac{1-r_w}{A_w}+\frac{4r_w}{9}\right)S.
\end{align*}
Together with \(\operatorname{vol}(P)\le C_wS\) from \cref{q:volume}, this proves the volume bound.

If the volume bound is attained, both inequalities in the preceding display must be equalities. The first forces
\[
\operatorname{area}(P_j)=A_w(k_j+1)
\qquad(j=1,w-1).
\]
Since \(\frac{4}{9}-\frac{1}{A_w}>0\), the second forces equality in \cref{g:share}, giving
\[
\frac{f(1)}{f(w-1)}\in\left\{w-1,\frac1{w-1}\right\}.
\]
The remaining assertion follows upon squaring this ratio.
\end{proof}

\begin{proposition}\label{g:combined}
Let $P$ be a lattice three-polytope with $I(P)=k$ and
$w=\lw(P)$. If $4\le w\le8$ and $k\ge w$, then
\[
 \vol(P)<6(k+1).
\]
\end{proposition}
\begin{proof}
By an affine unimodular change of coordinates, we may assume that the height of $P$ ranges over $[0,w]$, as in \eqref{eq:normalization}. If $P_j$ is contained in a split for some integer $1\le j\le w-1$, by \cref{pr:projection} and
\cref{pr:volume}, we get 
\[
 \vol(P)\le\frac{w^2}{(w-1)(w-2)}k+\frac{w^2}{2}
 \le\frac83k+4w<6(k+1).
\]
Indeed, the second inequality uses $4\le w\le8$, and the last
follows from
\[
 6(k+1)-\left(\frac83k+4w\right)
 =\frac{10}{3}(k-w)+6-\frac{2w}{3}>0.
\]

Suppose now that none of \(P_1,\ldots,P_{w-1}\) is contained in a split. By \cref{g:nonsplit}, we have 
\[
\operatorname{vol}(P)\le B_w(k+w-1)\le6(k+1).
\]
The second inequality follows by substituting \(w=4,5,6,7,8\) into the formula for \(B_w\) and using \(k\ge w\); it is strict unless \(w=k=4\).

It remains to rule out \(\operatorname{vol}(P)=6(k+1)\) when \(w=k=4\). Suppose, for a contradiction, that this occurs. The bound in \cref{g:nonsplit} is then attained. Thus its equality conditions give
\[
\frac{\operatorname{area}(P_1)}{\operatorname{area}(P_3)}
\in\left\{9,\frac{1}{9}\right\},
\qquad
\operatorname{area}(P_j)=A_4(k_j+1)
\quad(j=1,3).
\]
Consequently, $\frac{k_1+1}{k_3+1}\in\left\{9,\frac{1}{9}\right\}$, which implies \(\max\{k_1,k_3\}\ge8\), contradicting \(k_1+k_3\le k=4\). This completes the proof.
\end{proof}

\subsection{Proof of the main theorem}\label{sec:completion}

By \cref{g:combined}, we have \(\operatorname{vol}(P)<6(k+1)\) whenever \(4\le w\le8\) and \(k\ge w\). Moreover, \cref{s:parity} excludes \(w=8\) when \(k<8\). The remaining cases are therefore
\[
4\le w\le7,\qquad 3\le k<w.
\]
These are treated in \cref{s:section}, where each case is shown either to satisfy the strict inequality or to be impossible.

\begin{proof}[Proof of \cref{main:theorem}]
Since \(k\ge3\), the lattice width \(w\) is an integer at least two. For \(w=2\), by \cref{w2:theorem}, we get the volume bound and its equality characterization. For \(w=3\), the inequality is strict by \cref{w3:strict}. The cases \(4\le w\le8\) are covered by \cref{g:combined,s:parity} and \cref{s:section}. For \(w\ge9\), it suffices to use \cref{l:large}. Thus the inequality is strict whenever \(w\ge3\). Finally, \(S_{3,k}\) has exactly \(k\) interior lattice points and volume \(6(k+1)\), as verified in \cref{w2:theorem}.
\end{proof}

To resolve \cref{ConjectureZPW}, we first completed the proof of \cref{main:theorem} (i.e., the 3-dimensional lattice polytope case) and then attempted to settle the lattice simplex case. In \cite{LiuTangTaoZhang-Simplex}, we provided a complete proof of the lattice simplex case.

\appendix

\section{Proof of the sharp sampling inequality}
\label{q:proofsection}
Fix an integer $w\ge2$, and let
$C_w=\frac{2w^2}{(w-1)(2w-1)}$.
For bounded measurable functions $u,v$ on $[0,w]$, define
\[
\mathcal D(u,v)
=C_w\sum_{j=1}^{w-1}u(j)v(j)-\int_0^w u(t)v(t)\,dt.
\]
This expression is symmetric and linear in each argument.
For $0<a<w$, set
\[
G_a(t)=
\begin{cases}
t(w-a),&0\le t\le a,\\
a(w-t),&a\le t\le w.
\end{cases}
\]

\begin{lemma}\label{q:kernel}
For $0<a,b<w$, we have $\mathcal D(G_a,G_b)\ge0$.
\end{lemma}
\begin{proof}

Since $\mathcal D$ is symmetric, we may assume that
the original parameters satisfy $0<\alpha\le\beta<w$.
Setting $a=\alpha$ and $b=w-\beta$, we have
$a,b>0$ and $a+b=w-(\beta-\alpha)\le w$.
Thus it suffices to prove
$\mathcal D(G_a,G_{w-b})\ge0$ under these conditions.
 Set
\[
m=\lceil a\rceil-1,\qquad n=\lceil b\rceil-1.
\]
Then
\begin{equation}\label{q:index-bound}
m<a\le m+1,\qquad n<b\le n+1,\qquad m+n\le w-1.
\end{equation}
Integrating over the three intervals determined by the points \(0\), \(a\), \(w-b\), and \(w\), we obtain
\[
\frac1{ab}\int_0^wG_a(t)G_{w-b}(t)\,dt=\frac w6(w^2-a^2-b^2).
\]
For $a\le j\le w-b$, we have
\[
\frac{G_a(j)G_{w-b}(j)}{ab}=j(w-j).
\]
For $j<a$, the difference from $j(w-j)$ is
$w\left(\frac{j^2}{a}-j\right)$. Similarly, for $j>w-b$, the
difference is $w\left(\frac{(w-j)^2}{b}-(w-j)\right)$.
Since $m=\lceil a\rceil-1$ and $n=\lceil b\rceil-1$,
by summing over the three ranges, we get
\begin{align*}
\frac1{ab}\sum_{j=1}^{w-1}G_a(j)G_{w-b}(j)
=\frac{w(w^2-1)}6+w\sum_{j=1}^{m}\left(\frac{j^2}{a}-j\right)
+w\sum_{j=1}^{n}\left(\frac{j^2}{b}-j\right).
\end{align*}

Define
\[
F(t)=\frac{t^2}{6}
+C_w\sum_{j=1}^{\lceil t\rceil-1}
\left(\frac{j^2}{t}-j\right),\qquad 0<t\le w,
\]
where an empty sum is zero. We substitute the preceding expressions into the definition of $\mathcal D$ to obtain
\begin{equation}\label{q:kernel-defect}
\frac{\mathcal D(G_a,G_{w-b})}{wab}
=\frac{w^2}{2(2w-1)}+F(a)+F(b).
\end{equation}
For $1\le r\le w-1$ and $r<t\le r+1$, we have 
\[
F(t)=\frac{t^2}{6}+\frac{C_wr(r+1)(2r+1)}{6t}-\frac{C_wr(r+1)}2.
\]
By applying the elementary inequality
\[
t^2+\frac{A}{t}\ge3\left(\frac{A}{2}\right)^{2/3}
\qquad(A,t>0),
\]
we obtain
\begin{align*}
\frac{F(t)}{r(r+1)}\ge\frac{C_w}{2}
\left(
\left(\frac{1+\frac1{4r(r+1)}}{C_w}\right)^{1/3}-1
\right)\ge-\frac{w}{2(w-1)(2w-1)}.
\end{align*}
For the second inequality, the expression in brackets decreases with
$r$ and equals $-\frac{1}{2w}$ at $r=w-1$. Thus
\[
F(t)\ge-\frac{wr(r+1)}{2(w-1)(2w-1)}
\qquad(1\le r\le w-1,\ r<t\le r+1).
\]
This last bound also holds for $r=0$, since $F(t)=\frac{t^2}{6}\ge0$ when $0<t\le1$.
By \eqref{q:index-bound}, we have 
\[
m(m+1)+n(n+1)\le(m+n)(m+n+1)\le w(w-1).
\]
Consequently,
\[
F(a)+F(b)\ge-\frac{w\bigl(m(m+1)+n(n+1)\bigr)}{2(w-1)(2w-1)}\ge-\frac{w^2}{2(2w-1)}.
\]
The conclusion follows from \eqref{q:kernel-defect}.
\end{proof}

\begin{proof}[Proof of \cref{q:bilinear}]
We first extend \cref{q:kernel} to the functions $t$ and $w-t$.
For each fixed $0<t<w$, the definition of $G_a$ gives
\[
\frac{G_a(t)}a\longrightarrow w-t
\quad(a\downarrow0),\qquad
\frac{G_a(t)}{w-a}\longrightarrow t
\quad(a\uparrow w).
\]
Indeed, $\frac{G_a(t)}{a}=w-t$ on $[a,w]$, while
$\frac{G_a(t)}{w-a}=t$ on $[0,a]$. Both normalized functions are
bounded by $w$. For any fixed bounded second factor, the product
integrals therefore converge, since the intervals $[0,a]$ and
$[a,w]$ on which the functions differ shrink to zero length.
The interior integer samples eventually agree exactly with their limits.
We divide the inequality in \cref{q:kernel} by the appropriate positive factors and then take limits successively in the two arguments. This yields
\[
\mathcal D(u,v)\ge0
\]
for any $u,v$ chosen from $t$, $w-t$, and $G_a$ with $0<a<w$.

For a continuous nonnegative concave piecewise affine
function $f$, let $a_1,\ldots,a_N$ be the points where its
slope changes, and let $\delta_r\ge0$ be the decrease
in slope at $a_r$. Then
\[
f(t)=\frac{f(0)}w(w-t)+\frac{f(w)}w t
+\sum_{r=1}^{N}\frac{\delta_r}{w}G_{a_r}(t).
\]
Indeed, both sides have the same endpoint values and
changes in slope, since the slope of $G_a$ decreases
by $w$ at $a$. All coefficients are nonnegative. Thus the linearity in each argument proves the inequality for
piecewise affine $f,g$.

For continuous nonnegative concave $f,g$, join their
values on increasingly fine partitions by line segments.
These interpolants remain nonnegative and concave and
converge uniformly. Therefore, their product integrals and sampled
products converge. This proves the continuous case.

Finally, a finite nonnegative concave function on $[0,w]$
is continuous in $(0,w)$ and has one-sided limits at both
endpoints. The endpoint limits are finite, since the concavity and nonnegativity of $f$ imply
\[
0\le f(t)\le f(t)+f(w-t)\le2f\left(\frac w2\right).
\]
Moreover, the concavity of $f$ ensures that these limits are at least the corresponding endpoint values. We replace the endpoint values of $f$ by their limits to obtain a continuous nonnegative concave function with the same integral and values at $1,\ldots,w-1$. The same modification applies to $g$. This completes the proof for general finite nonnegative concave functions.

For $f(t)=g(t)=t$, we have 
\[
\frac{\int_0^w t^2\,dt}{\sum_{j=1}^{w-1}j^2}
=\frac{2w^2}{(w-1)(2w-1)}=C_w.
\]
Thus the constant is optimal.
\end{proof}

\section{Small interior-point cases}\label{s:section}

This appendix treats the remaining cases $4\le w\le7$ and $3\le k<w$. We use the section estimates above and, at widths $6$ and $7$, the additional results specified in the corresponding proofs.
Throughout this appendix, $P$ has the minimum-width normalization \eqref{eq:normalization}, $k_j=I(P_j)$, and
\[
 f(t)=\sqrt{\area(P_t)}\qquad(0\le t\le w).
\]
The function $f$ is continuous, nonnegative, and concave. By \cref{p:sections} and \cref{p:width}, we have
\begin{equation}\label{s:middle}
 \area(P_j)\le
 \begin{cases}
 2,&k_j=0,\\
 2k_j+\frac{5}{2},&k_j\ge1,
 \end{cases}
 \qquad 2\le j\le w-2.
\end{equation}
For $P_1$ and $P_{w-1}$, by \cref{p:rational}, we have 
\begin{equation}\label{s:outer}
 \area(P_j)\le\frac{(w+1)^2}{2w}(k_j+1),\qquad j=1,w-1.
\end{equation}
This holds provided that the section is not contained in a split.

At widths $4,5,6$, we may assume that no internal integer-height section is contained in a split. 
Otherwise, by \cref{pr:projection,pr:volume}, we obtain the following bounds for $w=4,5,6$, respectively,
\begin{equation}\label{s:split}
 \vol(P)\le\frac83(k+3),\qquad
 \vol(P)\le\frac{25}{12}(k+6),\qquad
 \vol(P)\le\frac95(k+10).
\end{equation}
The differences between $6(k+1)$ and these bounds are $\frac{10k-6}{3}$, $\frac{47k-78}{12}$, and $\frac{21k-60}{5}$, respectively.
They are positive for $k\ge3$.

The concavity and nonnegativity of $f$ imply
\begin{equation}\label{s:transfer}
 f(j)\le\frac{j}{i}f(i)\quad (0<i\le j<w),\qquad
 f(j)\le\frac{w-j}{w-i}f(i)\quad (0<j\le i<w).
\end{equation}
We shall also use the following elementary observation:
\begin{equation}\label{s:linear-product}
 px\le y\le qx,\quad 0<p\le q
 \quad\Longrightarrow\quad
 xy\ge\frac{pq\,x^2+y^2}{p+q},
\end{equation}
which follows by expanding $(qx-y)(y-px)\ge0$.
For adjacent internal samples $x=f(j)$, $y=f(j+1)$ with $1\le j\le w-2$, by \eqref{s:transfer}, we have 
$p=\frac{w-j-1}{w-j}$ and $q=\frac{j+1}{j}$.

\subsection{Width four}

\begin{proposition}\label{s:four}
If $\lw(P)=4$ and $k=3$, then $\vol(P)<24$.
\end{proposition}
\begin{proof}
By \eqref{s:split}, we may assume that none of $P_1,P_2,P_3$ is contained in a split. 
Therefore, we consider using \eqref{s:outer}. Set $a=f(1)$, $b=f(2)$ and $c=f(3)$.
If $k_2\le1$, then \cref{rem:midareatovol} and \eqref{s:middle} give
\[
 \vol(P)\le\frac{16}{3}\area(P_2)\le24.
\]
Suppose equality holds. Then $k_2=1$,
$\area(P_2)=\frac{9}{2}$, and the equality condition
in \cref{rem:midareatovol} forces $f$ to be
affine and vanish at one endpoint. Hence $\max\{\area(P_1),\area(P_3)\}=\frac{81}{8}$.
However, $k_1+k_3=2$ and \eqref{s:outer} imply $\max\{\area(P_1),\area(P_3)\}
\le\frac{75}{8}$, a contradiction. Thus $\vol(P)<24$ in this case.

Now suppose that $k_2\ge2$. By \cref{p:slabs}, we get
\begin{align*}
 3\vol(P)=&\area(P_0)+\area(P_4)+2(a^2+b^2+c^2)+(V_0+V_1)+(V_2+V_3)\\
 \le&(2a-b)^2+(2c-b)^2+2(a^2+b^2+c^2)+2a^2+2c^2\\
 =&8a^2+8c^2+4b^2-4b(a+c).
\end{align*}
The bounds in \eqref{s:transfer} give that $\frac{2a}{3}\le b\le2a$ and
$\frac{2c}{3}\le b\le2c$. Thus \eqref{s:linear-product} yields
$b(a+c)\ge\frac{a^2+c^2}{2}+\frac{3b^2}{4}$, and consequently
\begin{align*}
 3\vol(P)\le6(a^2+c^2)+b^2\le\frac{75}{4}(5-k_2)+2k_2+\frac52
 =\frac{385-67k_2}{4}\le\frac{251}{4}<72.
\end{align*}
Here we used $k_1+k_3=3-k_2$, \eqref{s:outer} for $a^2+c^2$, \eqref{s:middle} for $b^2$, and $k_2\ge2$.
This completes the proof. 
\end{proof}

\subsection{Width five}

\begin{lemma}\label{s:twonodes}
Let $P$ satisfy the hypotheses of \cref{p:slabs} with $w=5$. Let $f(t)=\sqrt{\area(P_t)}$. Then
\[
 \vol(P)\le\frac{35\bigl(f(2)^2+f(3)^2\bigr)-55f(2)f(3)}3.
\]
In particular, if $f(2),f(3)\le M$ for some $M\ge0$, then $\vol(P)\le5M^2$.
\end{lemma}
\begin{proof}
For convenience, we set $a=f(1)$, $b=f(2)$, $c=f(3)$, $d=f(4)$. By \cref{p:slabs}, we have
\[
V_0+V_1\le2a^2,\qquad V_3+V_4\le2d^2,
\]
and
\[
2V_2=(V_1+V_2)+(V_2+V_3)-V_1-V_3
\le2(b^2+c^2)-ab-cd.
\]
Hence
\[
\sum_{j=0}^4V_j
=(V_0+V_1)+V_2+(V_3+V_4)
\le2(a^2+d^2)+(b^2+c^2)-\frac12(ab+cd).
\]
The endpoint bounds also give
\[
\area(P_0)+\area(P_5)
\le(2a-b)^2+(2d-c)^2=4(a^2+d^2)+(b^2+c^2)-4(ab+cd).
\]
By substituting these estimates into \eqref{p:slabidentity}, we obtain
\begin{equation}\label{s:fivemixed}
3\vol(P)\le8(a^2+d^2)+4(b^2+c^2)-\frac92(ab+cd).
\end{equation}
The concavity and nonnegativity of $f$ further imply
\[
 \frac{b}{2}\le a\le2b-c,\qquad \frac{c}{2}\le d\le2c-b.
\]
For fixed $b,c$, the right-hand side of \eqref{s:fivemixed}
increases with $a,d$ on these intervals: its derivatives are
$16a-\frac{9b}{2}\ge\frac{7b}{2}>0$ and $16d-\frac{9c}{2}\ge\frac{7c}{2}>0$.
The upper bounds for $a$ and $d$ therefore yield
\[
3\vol(P)\le35(b^2+c^2)-55bc.
\]
Finally, suppose $b,c\le M$ and reflect height if necessary so that
$b\le c$. The concavity of $f$ implies $\frac{2c}{3}\le b\le c$. By
\eqref{s:linear-product} we obtain $bc\ge\frac{3b^2+2c^2}{5}$.
Consequently,
\[
3\vol(P)\le35(b^2+c^2)-55bc
\le2b^2+13c^2\le15c^2\le15M^2.
\]
This completes the proof.
\end{proof}

\begin{proposition}\label{s:five}
If $\lw(P)=5$ and $k\in\{3,4\}$, then $\vol(P)<6(k+1)$.
\end{proposition}
\begin{proof}
By \eqref{s:split}, we may assume that none of $P_1,\ldots,P_4$ is contained in a split, allowing us to apply \eqref{s:outer}. We adopt the notation $a,b,c,d$ from the proof of \cref{s:twonodes}.

If $\max\{b^2,c^2\}\le\frac92$, then \cref{s:twonodes} yields
\[
\vol(P)\le\frac{45}{2}<6(k+1)
\qquad(k\in\{3,4\}).
\]

Otherwise, we may assume that $b^2>\frac92$ by reversing the height coordinate if necessary. By \eqref{s:middle}, this implies $k_2\ge2$. Moreover, the concavity of $f$ gives
\[
c^2\ge\frac{4b^2}{9}>2,
\]
so \eqref{s:middle} also yields $k_3\ge1$. Consequently, $t=k_1+k_4\le k-3$.

The concavity of $f$ further implies
\[
\frac{3a}{4}\le b\le2a,\qquad
\frac{3d}{4}\le c\le2d.
\]
These bounds, together with \eqref{s:linear-product}, yield
\[
ab+cd\ge\frac{6(a^2+d^2)+4(b^2+c^2)}{11}.
\]
We now combine this estimate with \eqref{s:fivemixed} and the area bounds \eqref{s:outer} and \eqref{s:middle} to obtain

\begin{align*}
 3\vol(P)\le\frac{61}{11}(a^2+d^2)+\frac{26}{11}(b^2+c^2)
 &\le\frac{61}{11}\frac{18}{5}(t+2)+\frac{26}{11}\bigl(2(k-t)+5\bigr)
 \\& \le\frac{1098k+332}{55}<18(k+1).
 \end{align*}
The last two inequalities follow from $t\le k-3$ and $658-108k>0$ for $k\in\{3,4\}$, respectively.
\end{proof}

\subsection{Width six}

\begin{proposition}\label{s:six}
If $\lw(P)=6$ and $3\le k\le5$, then $\vol(P)<6(k+1)$.
\end{proposition}

\begin{proof}
By \eqref{s:split}, we may assume that no internal
integer-height section is contained in a split.
Let $A_j=\area(P_j)$ and set
\[
P^-=P\cap\{0\le z\le3\},\qquad
P^+=P\cap\{3\le z\le6\}.
\]
Suppose first that $k_3\ge2$.
By \eqref{s:transfer}, we have $\frac{4f(1)}{5}\le f(2)\le2f(1)$ and $\frac{4f(5)}{5}\le f(4)\le2f(5)$. Thus \eqref{s:linear-product} gives
\[
f(1)f(2)+f(4)f(5)\ge\frac47(A_1+A_5)+\frac5{14}(A_2+A_4).
\]
Both pieces have integer vertex heights and height span three. By applying \cref{w3:geometry} to $P^-$ and $P^+$,
we obtain
\[
\vol(P)\le
3\bigl((A_1+A_5)+A_2+A_4-f(1)f(2)-f(4)f(5)\bigr)\le\frac97(A_1+A_5)+\frac{27}{14}(A_2+A_4).
\]
According to the area bounds \eqref{s:outer} and \eqref{s:middle}, it follows that $A_1+A_5\le\frac{49(k_1+k_5+2)}{12}$ and $A_2+A_4\le\frac{5(k_2+k_4)}{2}+4$.
Therefore,
\begin{align*}
\vol(P)\le\frac{21}{4}(k_1+k_5+2)+\frac{135}{28}(k_2+k_4)
+\frac{54}{7}\le\frac{21}{4}k+\frac{54}{7}<6(k+1),
\end{align*}
where the last two inequalities use $k_1+k_5+k_2+k_4\le k-2$
and $k\ge3$, respectively.

Suppose now that $k_3\le1$. By \cref{p:sections}, we have $\lw(P_3)\ge3$, so \cite[Theorem~1.1]{ACFH2026} implies $P_3\simeq3\Delta_2$. We may apply a planar affine unimodular transformation while leaving the height coordinate unchanged, and thus assume that $P_3=3\Delta_2$. Consequently, $k_3=1$.

As $P_3$ is a lattice polygon, all new vertices created by the cut are lattice points. Hence both $P^-$ and $P^+$ are lattice polytopes. Furthermore, both have lattice width exactly $3$. Indeed, their common section $3\Delta_2$ has width at least $3$ in every integer direction with a nonzero horizontal component, while the height direction gives width $3$.

We claim that
\[
\vol(R)<6(I(R)+1)
\qquad(R=P^-,P^+).
\]
For $I(R)\ge3$, this follows from
\cref{w3:strict}. For $I(R)=1,2$, it follows
from the sharp volume bounds in \cite{Kas10} and
\cite[Corollary~5.2]{BK2016}: all equality polytopes have lattice
width two (see also \cite[Section~1.1]{BK2016}).
If $I(R)=0$, we may assume $R=P^-$ by reversing the height coordinate if necessary. Then $k_2=0$. By \eqref{s:middle} and the concavity of $f$, we obtain
\[
\sqrt2\ge f(2)\ge\frac13f(0)+\frac23f(3)=\frac13f(0)+\sqrt2.
\]
Equality forces $f(t)=\frac{t}{\sqrt2}$ on $[0,3]$,
and hence $\vol(P^-)=\frac{9}{2}<6$.
Thus the claim holds for both pieces.
Finally, since $k_3=1$, we have
$I(P^-)+I(P^+)=k-1$. Therefore, we have 
\begin{align*}
\vol(P)=\vol(P^-)+\vol(P^+)<6\bigl(I(P^-)+I(P^+)+2\bigr)=6(k+1).
\end{align*}
This completes the proof.
\end{proof}

\subsection{Width seven}
\begin{proposition}\label{s:seven}
There is no lattice three-polytope of width seven with exactly
three interior lattice points. If $\lw(P)=7$ and $4\le k\le6$,
then $\vol(P)<6(k+1)$.
\end{proposition}

\begin{proof}
We first show that
\begin{equation}\label{s:seven-central-budget}
k_3+k_4\ge4.
\end{equation}
Suppose instead that $k_3+k_4\le3$. By reversing the height coordinate if necessary, we may assume $k_3\le1$.
By \cref{p:sections} and \cite[Theorem~1.1]{ACFH2026},
a planar affine unimodular transformation allows us to take $P_3=3\Delta_2$. Thus $k_3=1$ and $k_4\le2$.
By \cite[Theorem~B]{ACFH2026}, there is a nonzero
integer direction $v$ with
$W_4:=\operatorname{width}_v(P_4)\le\frac{10}{3}$.

Set $W_t=\operatorname{width}_v(P_t)$. By the convexity of $P$ and the monotonicity and Minkowski additivity of width, $W_t$ is nonnegative and concave in $t$. Consequently,
\[
W_4\ge\frac34W_3+\frac14W_7\ge\frac34W_3,
\qquad
W_3\le\frac43W_4\le\frac{40}{9}<6.
\]
As $P_3=3\Delta_2$, we have $W_3\in3\Z_{>0}$,
and hence $W_3=3$. However, by \eqref{eq:interval} and \cref{p:sections}, we have
\[
7\le2W_{7/2}=W_3+W_4
\le3+\frac{10}{3}<7,
\]
a contradiction. Thus \eqref{s:seven-central-budget} holds, excluding $k=3$.

Now let $4\le k\le6$ and write $A_j=\area(P_j)$.
Cut $P$ at height three into $P^-$ and $P^+$.
By \eqref{s:transfer}, we have
\[
\frac12\le x:=\frac{f(1)}{f(2)}\le\frac65.
\]
The lower piece has height span three and integer vertex heights. Thus, by \cref{w3:geometry}, we have
\[
\vol(P^-)\le3A_2(x^2-x+1)
\le\frac{93}{25}A_2<4A_2,
\]
since $x^2-x+1$ is nondecreasing for $x\ge\frac{1}{2}$.
The midpoint estimate \eqref{w2:rescaled-middle}
on $P^+$ gives $\vol(P^+)\le\frac{16A_5}{3}$.
Thus
\[
\vol(P)<4A_2+\frac{16}{3}A_5.
\]
The reflection of the height coordinate yields the same bound with $A_2$ and $A_5$ interchanged. By taking the average of these two bounds and using $k_2+k_5\le k-4$ and \eqref{s:middle}, we obtain
\[
\vol(P)<\frac{14}{3}(A_2+A_5)\le\frac{14}{3}(2k-3)\le6(k+1),
\]
where the last inequality uses $k\le6$.
This completes the proof.
\end{proof}






\vspace{5mm}
\noindent
{\small \textbf{Acknowledgments:}}
Feihu Liu was partially supported by the Postdoctoral Fellowship Program and China Postdoctoral Science Foundation (Grant No. BX2026002).

\noindent{\small \textbf{Declaration of AI Assistance:}}
During the preparation of this manuscript, the authors used ChatGPT to assist in exploring possible approaches, checking technical details, and improving the exposition. All mathematical arguments, computations, proofs, and results were independently verified by the authors. The authors take full responsibility for the content of this manuscript.

\end{document}